\documentclass[a4paper,12pt,reqno]{amsart}
\usepackage{amsfonts}
\usepackage{amsmath}
\usepackage{amssymb}
\usepackage{mathrsfs}
\usepackage{hyperref}
\usepackage[english,francais]{babel}
\usepackage{multicol}

\newtheorem{e-proposition}[theorem]{Proposition}

\newtheorem{e-definition}[theorem]{Definition\rm}

 \newtheorem{thm}{Theorem}[section]
 \newtheorem{cor}[thm]{Corollary}

 \theoremstyle{definition}
 \newtheorem{defn}[thm]{Definition}
\newcommand{\er}{\mathbb{R}}
\newcommand{\M}{\mathcal{M}}

\newcommand{\zn}{\mathbb{Z}^n}
\newcommand{\bba}{\mathcal{B}}
\newcommand{\W}{\mathcal{W}}
\newcommand{\Rev}{\mathcal{R}}
\newcommand{\Tr}{Tr}
\newcommand{\tn}{\mathbb{T}^n}
\def\Rn{{{\mathbb R}^n}}
\def\Tn{{{\mathbb T}^n}}
\def\Zn{{{\mathbb Z}^n}}

\newcommand{\ern}{{\mathbb{R}}^n}

\newcommand{\bi}{\begin{itemize}}
\newcommand{\ei}{\end{itemize}}
\newcommand{\be}{\begin{enumerate}}
\newcommand{\ee}{\end{enumerate}}
\newcommand{\ard}{{\mathbb{R}}^{2d}}
\newcommand{\arda}{{\mathbb{R}}^{d}}
\newcommand{\beq}{\begin{equation}}
\newcommand{\eq}{\end{equation}}
 
\newcommand{\efet}{{\mathcal{F}}_{\mathbb{T}^n}}

\newcommand{\vel}{L^{p(\cdot)}}
\newcommand{\velp}{L^{p'(\cdot)}}
\newcommand{\velq}{L^{q(\cdot)}}

\renewcommand{\Tr}{{\rm Tr}}
\begin{document}
\selectlanguage{english}
\title[Grothendieck-Lidskii formulae]
{Grothendieck-Lidskii trace formula for mixed-norm and variable Lebesgue spaces}

\author[Julio Delgado]{Julio Delgado}


\address{Department of Mathematics\\
Imperial College London\\
180 Queen's Gate, London SW7 2AZ\\
United Kingdom}

\email{j.delgado@imperial.ac.uk}

\author{Michael Ruzhansky}

\address{Department of Mathematics\\
Imperial College London\\
180 Queen's Gate, London SW7 2AZ\\
United Kingdom}

\email{m.ruzhansky@imperial.ac.uk}

\author{Baoxiang Wang}

\address{LMAM, School of Mathematical Sciences\\
Peking University\\
Beijing 100871\\
China}
\email{wbx@pku.edu.cn}



\date{\today}

\thanks{The first author was supported by the
Leverhulme Research Grant RPG-2014-02.
The second author was supported by
 the EPSRC Grant EP/K039407/1. No new data was collected or generated during the course of the research.
 The authors have been also supported by the Sino-UK research project by the British
Council China and the China Scholarship Council, and by EPSRC Mathematics Platform grant EP/I019111/1.}

\begin{abstract}
In this note we present the metric approximation property for weighted mixed-norm 
$L_w^{(p_1,\dots ,p_n)}$ and variable exponent Lebesgue type spaces. 
As a consequence, this also implies the same property for modulation and Wiener-Amalgam spaces.
We then characterise nuclear operators on such spaces and state the corresponding Grothendieck-Lidskii trace formulae.
 We apply the obtained results to derive criteria for nuclearity and trace formulae
for periodic operators on $\Rn$ and functions of the harmonic oscillator in terms of global symbols.
\end{abstract}



\dedicatory{To the memory of Yuri Safarov}

\maketitle

\section{Introduction}

The famous Lidskii formula \cite{li:formula} states the equality between the operator trace of trace-class operators in Hilbert spaces and the sum of their eigenvalues. Similar properties can be established in Banach spaces using Grothendieck's theory of nuclear operators. The basic ingredient for such an approach is the approximation property of the Banach space under consideration. In this paper we describe the approximation properties and the subsequent Grothendieck-Lidskii formulae for the traces of operators in weighted mixed-norm and in variable exponent Lebesgue spaces.

The approximation property is one of the fundamental properties of the `geometry' of Banach spaces. A particular importance of this property from the Grothendieck's viewpoint is that once a Banach space is known to have it, the trace can be defined and consequently the Fredholm's determinant leading to numerous further developments. Thus, the topic finds itself closely related to a wide range of analysis: spectral analysis, operator theory, functional analysis, harmonic analysis, PDEs. In particular, nuclearity properties of a given operator allow to obtain information on the distribution eigenvalues
 through the use of Grothendieck, K\"onig, Maurey, Retherford and Johnson inequalities which can be seen as an extension of the classical Weyl inequality relating sums of eigenvalues and Schatten-von Neumann norms.

In the paper we present the metric approximation property for three scales of spaces that are of importance in a broad range of mathematical subjects.
First, the mixed Lebesgue spaces are a basic tool for harmonic analysis and evolutions PDEs (Strichartz estimates). The approximation property of such spaces (\cite{drwap2:app}) gives rise to an introduction of spectral methods (following Grothendieck) to a variety of questions of harmonic analysis and PDEs. A part of the paper is also devoted to the development of these ideas. Second, Wiener amalgam spaces are a central object of the time-frequency analysis, another area with links to several mathematical subjects as well as its applications. Also, the approximation property in the scale of modulation spaces gives rise to the introduction of spectral analysis to PDEs of very different type - these spaces become more and more effective (compared to Besov spaces) in many types of PDEs including such PDEs as the Navier-Stokes equation, e.g. \cite{it:nst}. Finally, we discuss (following \cite{drap1:app})
 the validity of the metric approximation property for the variable exponent Lebesgue spaces. In all the considered spaces a characterisation of nuclear operators is established and the corresponding Grothendieck-Lidskii formulae for the trace are derived. In the case of variable exponent Lebesgue spaces on the torus we present sufficient conditions to ensure the nuclearity of an operator of the form $\alpha(x) (I-\Delta)^{-\frac{\tau}{2}}$,  where $\alpha$ is a suitable function.
  
Let $\bba_1, \bba_2$ be Banach spaces and let $0<r\leq 1$. A linear operator $T$
from $\bba_1$ into $\bba_2$ is called {\em r-nuclear} if there exist sequences
$(x_{n}^{\prime})\mbox{ in } \bba_1' $ and $(y_n) \mbox{ in } \bba_2$ so that
$
Tx= \sum\limits_{n=1}^{\infty} \left <x,x_{n}'\right>y_n \,\mbox{ and }\,
\sum\limits_{n=1}^{\infty} \|x_{n}'\|^{r}_{\bba_1'}\|y_n\|^{r}_{\bba_2} < \infty.
$
When $r=1$ the notion of nuclear operators agrees with the one of trace class operators
 in the setting of Hilbert spaces ($\bba_1=\bba_2=H$). Grothendieck \cite{gro:me} proved that the trace $\Tr(T)$ is well defined for 
 all nuclear operators if and only if the Banach space
 $\bba$ has the {\em approximation property},  i.e. for every compact
   set $K$ in $\bba$ and for every $\epsilon >0$, there exists $F\in \mathcal{F}(\bba) $ such that
$\|x-Fx\|<\epsilon ,\quad \textrm{ for all } x\in K,$
where $\mathcal{F}(\bba)$ denotes the space of finite rank bounded linear operators
 on $\bba$. If in the definition above
 the operator $F$ satisfies $\|F\|\leq 1$ one says that $\bba$ possesses the 
 {\em metric approximation property}. In \cite{gro:me} Grothendieck proved that if $T$ is $\frac 23$-nuclear from $\bba$ into $\bba$, then its trace agrees with the sum of eigenvalues. The first example of a Banach space without the approximation property appeared in \cite{ap:enflo}. In \cite{ap:sz} Szankowski proved that $B(H)$ does not have the approximation property. More recently, these properties have been intensively also investigated in different scales of function spaces, see e.g.
\cite{Alberti-Csornyei-Pelczynski-Preiss:BV}, \cite{jsza:happy},  
 \cite{lo:bapp1}.  
\section{Mixed-Normed $L^p$, Modulation and Wiener-Amalgam spaces}  
The modulation spaces have been intensively investigated in the last decades.
 Modulation spaces start finding numerous applications in various problems in linear
and nonlinear partial differential equations, see \cite{Ruzhansky-Sugimoto-Wang:B-2012}
for a recent survey. 
For a suitable weight $w$ on $\ard$, $1\leq p,q<\infty$ and a window $g\in\mathcal{S}(\arda)$ the modulation space $\M_{w}^{p,q}(\arda)$ consists of the temperate distributions $f\in\mathcal{S} '(\arda)$ such that
\begin{equation}\label{EQ:modul}
\|f\|_{\M_{w}^{p,q}}:=\|V_gf\|_{L_w^{p,q}}:=\end{equation}
\[\left(\int_{\arda}\left(\int_{\arda}|V_gf(x,\xi)|^pw(x,\xi)^pdx\right)^{\frac qp}
d\xi\right)^{\frac 1q}<\infty ,\]
where $V_gf(x,\xi)$ denotes the short-time Fourier transform of $f$ with respect to $g$ at the point $(x,\xi)$. The modulation space $\M_{w}^{p,q}(\arda)$ endowed with the above norm becomes a Banach space,
independent of $g\not=0$. 
A {\em weight function} is a non-negative, locally integrable function on $\ard$. 
%
A weight function $v$ on $\ard$ is called {\em submultiplicative, } if
\beq \label{polyw0} 
v(x+y)\leq v(x)v(y) \mbox{ for all }x,y\in\ard . 
\eq
A weight function $w$ on $\ard$ is  {\em v-moderate, } if
\beq 
w(x+y)\leq v(x)w(y) \mbox{ for all }x,y\in\ard . 
\eq
In particular the weights of polynomial type play an important role. They are of the form
\beq\label{polyw} 
v_s(x,\xi)=(1+|x|^2+|\xi|^2)^{s/2}.
\eq 
The $v_{s}$-moderated weights (for some $s$) are called {\em polynomially moderated}. 

We now recall the definition of weighted mixed-norm $L^p$ spaces. Let $(\Omega_i, S_i,\mu_i)$, for $i=1,\dots,n$, be given $\sigma$-finite measure spaces. We
write $x=(x_{1},\ldots,x_{n})$, and let
$P=(p_1,\dots,p_n)$ a given $n$-tuple with $1\leq p_i<\infty$. We say that $1\leq P<\infty $ if $1\leq p_i<\infty$ for all $i=1,\dots, n$. Let $w$ be a strictly positive measurable function. The  norm $\|\cdot\|_{L_w^P}$
 of a measurable function $f(x_1,\dots,x_n)$ on the corresponding product measure space is defined by
\[\|f\|_{L_w^P}:=\left(\int_{\Omega_n}\cdots\left(\int_{\Omega_2}\right.\right.\]
\[\left. \left.\left(\int_{\Omega_1}|f(x)|^{p_1}w(x)d\mu_1(x_{1})\right)^{\frac{p_2}{p_1}}d\mu_2(x_{2})\right)^{\frac{p_3}{p_2}}\cdots d\mu_n(x_{n})\right)^{\frac{1}{p_n}}.\]
$L_w^P$-spaces endowed with the $\|\cdot\|_{L_w^P}$-norm become Banach spaces and the dual $(L_w^P)'$ of $L_w^P$ is $L_{w^{-1}}^{P'}$, where $P'=(p_1',\dots,p_n')$. In view of our application to the modulation spaces $\M_w^{p,q}$ we will consider in particular the case of the index of the form $(P,Q)=(p_1,\dots,p_d,q_1,\dots,q_d)$ where $p_i=p, q_i=q$ and $\Omega_i=\er$ endowed with the Lebesgue measure. In this case the weight is taken in the form $w=w(x,\xi)$ where 
$x\in\er^d, \xi\in\er^d$. 
 In the rest of  this section we will assume that 
 the weights $w$ satisfy the following condition for all $x\in\Omega$:
 \beq\label{conw1}w(x_1,\dots,x_n)\leq w_1(x_1)\cdots w_n(x_n),\eq
where $w_j$ is a weight on $\Omega_j$
(i.e. a strictly positive locally integrable function).
In particular, the condition holds for polynomially moderate
weights on $\ern$ satisfying  for a suitable $n$-tuple $(\beta_1,\dots,\beta_n)$
the condition
\beq
\label{conw2}
w(x_1,\dots,x_n)\leq \langle x_1\rangle^{\beta_1}\cdots\langle x_n\rangle^{\beta_n},
\eq
where $\langle x_{j}\rangle =1+|x_{j}|$. The following theorem was established in
\cite{drwap2:app}.

\begin{thm} \label{app12a1} 
The weighted mixed-norm spaces $L_w^P=L_w^{(p_1,\dots ,p_n)}$  with $w$ satisfying (\ref{conw1}) have the metric approximation property.
\end{thm}
Let us recall now the definition of the Wiener amalgam spaces $\W^{p,q}_{w}(\mathbb R^{d})$.
There are several definitions possible for the spaces $\W^{p,q}_{w}$, in particular
involving the short-time Fourier transform similarly to the definition of the modulation
spaces in (\ref{EQ:modul}). 
To make an analogy with modulation spaces, we can reformulate their definition
 (\ref{EQ:modul}) in terms of the mixed-normed Lebesgue spaces, by saying that
\begin{equation}\label{EQ:def-mod}
f\in \M^{p,q}_{w}(\mathbb R^{d}) \;\textrm{ if and only if }\;
V_{g}f\cdot w\in L^{(p,q)}(\mathbb R^{d}\times \mathbb R^{d}).
\end{equation}
Now, for a function
$F\in L^{1}_{loc}(\mathbb R^{2d})$,
we denote $\Rev F(x,\xi):=F(\xi,x).$
Then we can define
\begin{equation}\label{EQ:def-was}
f\in \W^{p,q}_{w}(\mathbb R^{d}) \;\textrm{ if and only if }\;
\Rev (V_{g}f\cdot w)\in L^{(q,p)}(\mathbb R^{d}\times \mathbb R^{d}).
\end{equation}
However, for our purposes the following description
through the Fourier transform will be more practical. 
For a review of different definitions we refer to 
\cite{Ruzhansky-Sugimoto-Toft-Tomita:MN-2011}.
So, in what follows, we will
always assume that the weights in modulation and Wiener amalgam spaces
are submultiplicative and polynomially moderate (but we do not need to assume this when
talking about weighted mixed-norm $L^{P}$-spaces) as in
(\ref{polyw0})--(\ref{polyw}).
Then, because of the identity
$$
|V_{g}f(x,\xi)|=(2\pi)^{-d}\, |V_{\widehat{g}} \widehat{f}(\xi,-x)|,
$$
the Wiener amalgam space
$\W^{p,q}_{w}$ and the modulation spaces are related through the Fourier transform 
by the formula
\begin{equation}\label{EQ:WM}
\|f\|_{\W^{p,q}_{w}}\simeq \|\widehat{f}\|_{\M^{q,p}_{w_{0}}},
\end{equation}
 where $w(x,\xi)=w_{0}(\xi,-x).$
 As a consequence of Theorem \ref{app12a1} we now immediately obtain:
 
\begin{cor} \label{COR:mod}
Let $1\leq p,q<\infty$, and $w$ a submiltiplicative polynomially moderate weight. 
Then $\M_{w}^{p,q}$ has the metric approximation property. 
Consequently, also the Wiener amalgam space $\W^{p,q}_{w}$ 
has the metric approximation property.
\end{cor} 

It was observed by Feichtinger and Gr\"ochenig \cite{Feichtinger-Grochenig:approx-MM-1989} that the metric approximation property could be alternatively established for the corresponding sequence spaces using appropriate atomic decompositions. Thus, Corollary \ref{COR:mod} implies the metric approximation property for the sequence spaces arising through the atomic decompositions of $\M_{w}^{p,q}$ and $\W^{p,q}_{w}$. 

 In order to formulate a characterisation of $r$-nuclear operators between weighted mixed-norm spaces we will consider 
$1\leq P, Q<\infty$. The multi-index $P$ will be associated to 
 the measures $\mu_i$ $(i=1,\dots, l)$ and $Q$ will correspond to the measures 
$\nu_j$ $(j=1,\dots, m)$. We will also denote $\mu:=\mu_1\otimes\cdots\otimes\mu_l$ and 
 $\nu:=\nu_1\otimes\cdots\otimes\nu_m$ the corresponding product measures on the product spaces 
$\Omega=\prod\limits_{i=1}^l\Omega_i, \Xi=\prod\limits_{j=1}^m\Xi_j$. For a weight $w$ we will 
denote $w_P(\Omega):=\|1_{\Omega}\|_{L_w^P(\mu)}$. The additional property (\ref{conw1})
will be only required for the formulation of trace relations. 
\begin{defn}\label{triple} 
Let $({\Omega}_i,{\mathcal{M}}_i,\mu_i ) (i=1,\dots, l)$ be measure spaces and
$\mu:=\mu_1\otimes\cdots\otimes\mu_l$ the corresponding product measure on $\Omega=\prod\limits_{i=1}^l\Omega_i$
. We will also call $\Lambda\in\mathcal{M}:=\bigotimes\limits_{i=1}^l{\mathcal{M}}_i$ a {\em box} if it is of the form $\Lambda=\prod\limits_{i=1}^l\Lambda_i$. For a measure $\mu$, a weight $w$ on $\Omega$ and a multi-index $P$ we will say that the triple $(\mu,w, P)$ is $\sigma$-finite if there exists a family of disjoint boxes $\Omega^k$ such that $\mu(\Omega^k)<\infty$, $\bigcup\limits_{k=1}^{\infty}\Omega^k=\Omega$ and 
\[w_P(\Omega^k)=\|1_{\Omega^k}\|_{L_w^P(\mu)}<\infty.\]
\end{defn}

We can now give a characterisation of $r$-nuclear operators on weighted mixed-norm spaces and a trace formula. 
\begin{thm}\label{ch2} 
Let $0<r\leq 1$. Let $\,({\Omega}_i,{\mathcal{M}}_i,\mu_i ) \,\,(i=1,\dots, l)$, 
$({\Xi}_j,{\mathcal{M}}'_j,{\nu}_j) (\,j=1,\dots, m)$ be measure spaces. Let $1\leq P,Q<\infty$. 
Let $w, \widetilde{w}$ be weights on $\Omega, \Xi$ respectively such that the triples 
$(\mu,w, P)$, $(\nu,{\widetilde{w}}^{-1}, Q')$ are $\sigma$-finite. Then $T$ is
  $r$-nuclear operator from $L_w^{P}(\mu)$ into $L_{\widetilde{w}}^{Q}(\nu)$ if and only if there exist a sequence
 $(g_n)$ in $L_{\widetilde{w}}^{Q}(\nu)$, and a sequence $(h_n)$ in $L_{w^{-1}}^{P'}(\mu)$ such that $\sum \limits_{n=1}^\infty \|
 g_n\|_{L_{\widetilde{w}}^{Q}(\nu)}^r\|h_n\|_{L_{w^{-1}}^{P'}(\mu)}^r<\infty$, and such that for all $f\in L_w^{P}(\mu)$
\[Tf(x)=\int\limits_{\Omega}\left(\sum\limits_{n=1}^{\infty}
  g_n(x)h_n(y)\right)f(y)d\mu(y), \,\,\mbox{for a.e } x.\]
Moreover, if $w=\widetilde{w}$ satisfies (\ref{conw1}), $\mu=\nu$, $P=Q$ and $T$ is $r$-nuclear in $\mathcal{L}(L_w^{P}(\mu))$ with $r\leq \frac 23$, then  
\[\Tr(T)=\sum\limits_{j=1}^{\infty}\lambda_j,\]
where $\lambda_j\,\, (j=1,2,\dots)$ are the eigenvalues of $T$ with multiplicities taken into account,
and $\Tr(T)=\sum\limits_{j=1}^{\infty} \left< u_{j},v_{j} \right>.$
\end{thm}
\medskip
Analogous characterisations of $r$-nuclear operators can be obtained in the case of modulation spaces and Wiener-Amalgam spaces (cf. \cite{drwap2:app}). We now formulate an application in the case of modulation spaces to the study of functions of the harmonic oscillator $A=-\Delta +|x|^2$ on $\arda$, defined by 
\begin{equation}\label{EQ:HA-F}
F(-\Delta +|x|^2)\phi _j= F(\lambda_j)\phi _j, \qquad
j=1,2,\ldots,
\end{equation}
where $\lambda_j$'s are the eigenvalues of $A$.
We have:
\begin{thm} 
Let $0<r\leq 1$, $s\in\er$ and $1\leq p,q<\infty$. 
The operator $F(-\Delta +|x|^2)$ is $r$-nuclear on $\M_s^{p,q}(\arda)$ provided that
\begin{equation}\label{EQ:nucl}
\sum\limits_{j=1}^{\infty} |F(\lambda_j)|^r\|\phi _j
\|_{\M_s^{p,q}}^r\|\phi _j\|_{\M_{-s}^{p',q'}}^r<\infty.
\end{equation}
Moreover, if (\ref{EQ:nucl}) holds with $r=1$,
we have the trace formula
\begin{equation}\label{EQ:F-trace}
\Tr F(-\Delta +|x|^2) =\sum_{j=1}^{\infty} F(\lambda_{j}),
\end{equation}
with the absolutely convergent series.
\end{thm} 
\section{Variable exponent Lebesgue spaces}  
The variable exponent Lebesgue spaces are a generalisation of the classical Lebesgue spaces, replacing the constant exponent $p$ by a variable exponent function $p(x)$. The development of the analysis of many problems on those spaces has been of great interest in the last decades as has been exhibited in the recent books  \cite{di:book}, \cite{lpv:cufrw} and the literature therein. 
 We now recall the definition of variable exponent Lebesgue spaces and refer the reader to \cite{di:book} for the basic properties of such spaces. Let $(\Omega,\mathcal{M},\mu)$ be a $\sigma$-finite, complete measure space. We
define $\mathcal{P}(\Omega,\mu)$ to be the set of all $\mu$-measurable functions $p: \Omega\rightarrow [1,\infty].$ The functions in
$\mathcal{P}(\Omega,\mu)$ are called variable exponents on $\Omega$. We define $p^+=p_{\Omega}^+:={\rm ess}\sup_{x\in \Omega}p(x), \quad p^-=p_{\Omega}^-:={\rm ess}\inf_{x\in\Omega}p(x).$ If $p^+<\infty$, then $p$ is called a {\em bounded variable exponent}.
 If $f:\Omega\rightarrow\er$ is a measurable function we define the {\em modular} associated with $p=p(\cdot)$ by 
\[\rho_{p(\cdot)}(f):=\int\limits_{\Omega}|f(x)|^{p(x)}d\mu(x),
\]
and $\|f\|_{\vel(\mu)}:=\inf\{\lambda >0:\rho_{p(\cdot)}(f/\lambda)\leq 1\}.$ The resulting spaces $L^{p(\cdot)}(\mu)$ of measurable functions such that $\|f\|_{\vel(\mu)}<\infty$ are Banach spaces and enjoy many properties similar to the classical Lebesgue $L^p$ spaces. 
 If the variable exponent $p(\cdot)$ is bounded the space $L^{p(\cdot)}(\mu)$ is separable and if we denote by $p'(\cdot)$ the variable exponent defined pointwise by $\frac{1}{p(x)}+\frac{1}{p'(x)}=1,$
then $(L^{p(\cdot)}(\mu))'=L^{p'(\cdot)}(\mu)$. Moreover, if $1<p^-\leq p^+<\infty$ the space $L^{p(\cdot)}(\mu)$ is reflexive.  For the study of the approximation property we will restrict to consider bounded variable exponents due to the density of the simple functions in $L^{p(\cdot)}$ in that case.  
We can now state the metric approximation property which was proved in \cite{drap1:app}:

\begin{thm} \label{app12} Let $p\in\mathcal{P}(\Omega,\mu)$ be a bounded variable exponent. Then, the variable exponent Lebesgue space $L^{p(\cdot)}(\mu)$ has the metric approximation property.
\end{thm}

We are now ready to give a characterisation of $r$-nuclear operators for variable exponent spaces.
\begin{thm}\label{ch2}  
Let $({\Omega},{\mathcal{M}},\mu)$ and
$({\Xi},{\mathcal{M}}',\nu)$ be $\sigma$-finite complete measure spaces. Let $0<r\leq 1$. Then $T$ is
  $r$-nuclear operator from $\vel(\mu)$ into $\velq(\nu)$ if and only if there exist a sequence
 $(g_n)$ in $\velq(\nu)$, and a sequence $(h_n)$ in $\velp(\mu)$ such that $\sum \limits_{n=1}^\infty \|
 g_n\|_{\velq(\nu)}^r\|h_n\|_{\velp(\mu)}^r<\infty$, and such that for all $f\in\vel(\mu)$ we have
\[
Tf(x)=\int\limits_{\Omega}\left(\sum\limits_{n=1}^{\infty}
  g_n(x)h_n(y)\right)f(y)d\mu(y), \,\,\mbox{for a.e } x.
\]

Moreover, if $\Omega=\Xi$, $\mu=\nu$, $p(\cdot)=q(\cdot)$, $p^+<\infty$, and $T$ is $r$-nuclear in 
$\mathcal{L}(\vel(\mu))$ with $r\leq \frac 23$, then  
\[\Tr(T)=\sum\limits_{j=1}^{\infty}\lambda_j,\]
where $\lambda_j\,\, (j=1,2,\dots)$ are the eigenvalues of $T$ on $\vel(\mu)$ with multiplicities taken into account, 
and $$\Tr(T)=\sum\limits_{n=1}^{\infty} \left< g_n,h_n \right>=\int\limits_{\Omega}\sum\limits_{n=1}^{\infty}g_n(x)h_n(x)d\mu.$$
\end{thm}
We denote the $n$-dimensional torus by $\Tn=\Rn/ \Zn$. Its unitary dual can be 
described as $\widehat{\Tn}\simeq \Zn$, and the collection $\{ \xi_k(x)=e^{2\pi i x\cdot k}\}_{k\in\Zn}$ 
is an orthonormal basis of $L^2(\Tn)$. A corresponding operator is associated to a symbol $\sigma(x,\xi)$ which will be called a periodic pseudo-differential operator or the operator given by the toroidal quantization: 
\beq  T_{\sigma}f(x)=\sum\limits_{\xi\in\zn}
  e^{2\pi i x\cdot\xi}\sigma(x,\xi)(\efet f)(\xi) ,\eq
 which can also be written as 
\beq  \label{EQ:Tq}
T_{\sigma}f(x)=\sum\limits_{\xi\in\zn}\int\limits_{\tn}
  e^{2\pi i(x-y)\cdot\xi}\sigma(x,\xi)f(y)dy.
 \eq
We refer to \cite{Ruzhansky-Turunen-JFAA-torus} for an extensive analysis of such toroidal quantization, and to \cite{rt:book} for the toroidal background analysis. In the rest of this section we will consider $\tn$ endowed with the Borel $\sigma$-algebra and the Lebesgue measure so that we will just write $\mathcal{P}(\tn)$ 
 to denote the corresponding class of variable exponents. 
Given a measurable function $\alpha$ on $\tn$,
we take the symbols $\alpha(x)$ and  $\sigma(\xi)$, the corresponding multiplication is the operator denoted by 
$\alpha T_{\sigma}$ given by $\alpha T_{\sigma}f=\alpha\sigma(D)f$ on $\tn$. 
\begin{cor} 
Let $p(\cdot)\in\mathcal{P}(\tn)$. Let $0<r\leq 1$, $\alpha\in\velp$, and let $\sigma(\xi)$ be a symbol such that
\[\sum\limits_{\xi\in\zn}|\sigma(\xi)|^r<\infty.\]
Then $\alpha T_{\sigma}$ is $r$-nuclear from $\vel$ to $\velq$ for all $q(\cdot)\in\mathcal{P}(\tn)$. 
If additionally $p^+<\infty$,  $r\leq \frac 23$, and $q(\cdot)=p(\cdot)$, then 
$\alpha T_{\sigma}$ is $r$-nuclear in $\mathcal{L}(\vel(\tn))$ and   
\[
\Tr(\alpha T_{\sigma})=\int_{\tn}\alpha(x)dx\cdot\sum\limits_{\xi\in\zn} \sigma(\xi)=\sum\limits_{j=1}^{\infty}\lambda_j,
\]
where $\lambda_j\,\, (j=1,2,\dots)$ are the eigenvalues of $\alpha T_{\sigma}$ with multiplicities taken into account.
\end{cor}
\medskip

In particular, let us consider the symbol
$\sigma(\xi)=(1+4\pi^2 |\xi|^2)^{-\frac{\tau}{2}}$ for $\tau>0$.
The corresponding multiplication yields the operator $\alpha T_\sigma f=\alpha(I-\Delta)^{-\frac{\tau}{2}}f$ on $\tn$. 
We observe that $\sum\limits_{\xi\in\zn}(1+4\pi^2 |\xi|^2)^{-\frac{r\tau}{2}}<\infty$ if and only if $r\tau>n$. Consequently we obtain:

\begin{cor} 
Let $p(\cdot)\in\mathcal{P}(\tn)$. If $0<r\leq 1$, $\alpha\in\velp$, and $r\tau>n$, then 
 $\alpha T_{\sigma}=\alpha(I-\Delta)^{-\frac{\tau}{2}}$ is $r$-nuclear from $\vel$ to $\velq$ for all 
 $q(\cdot)\in\mathcal{P}(\tn)$. 
If additionally $p^+<\infty$,  $r\leq \frac 23$, and $q(\cdot)=p(\cdot)$, then $\alpha (I-\Delta)^{-\frac{\tau}{2}}$ is $r$-nuclear in 
$\mathcal{L}(\vel(\tn))$ and   
\[
\Tr(\alpha (I-\Delta)^{-\frac{\tau}{2}})=\int_{\tn}\alpha(x)dx\cdot\sum\limits_{\xi\in\zn}(1+4\pi^2|\xi|^2)^{-\frac{\tau}{2}}=
\sum\limits_{j=1}^{\infty}\lambda_j,\]
where $\lambda_j\,\, (j=1,2,\dots)$ are the eigenvalues of $\alpha (I-\Delta)^{-\frac{\tau}{2}}$ 
on $\vel(\tn)$ with multiplicities taken into account.
\end{cor} 

\end{document}